\documentclass[11pt]{35}
\usepackage{amsthm, amsmath, amssymb, amsbsy, amsfonts, latexsym, euscript}
\usepackage{calrsfs} 
\usepackage{graphicx}

\DeclareMathAlphabet{\varmathbb}{U}{pxsyb}{m}{n}

\def\leq{\leqslant}

\def\phi{\varphi}

\def\bar{\overline}

\def\kappa{\varkappa}

\newcommand{\D}{\mathrm{d}\kern0.2pt}%

\newcommand{\E}{\mathrm{e}\kern0.2pt} 
\newcommand{\ii}{\kern0.05em\mathrm{i}\kern0.05em}

\newcommand{\RR}{\mathbb{R}}%
\newcommand{\DD}{\varmathbb{D}}%

\theoremstyle{remark}

\numberwithin{equation}{section}

\begin{document}

\noindent {\Large \bf Sloshing in containers with vertical walls: isoperimetric
\\[3pt] inequalities for the fundamental eigenvalue}

\vskip5mm

{\bf Nikolay Kuznetsov}

\vskip-2pt {\small Laboratory for Mathematical Modelling of Wave Phenomena}
\vskip-4pt {\small Institute for Problems in Mechanical Engineering, RAS} \vskip-4pt
{\small V.O., Bol'shoy pr. 61, St Petersburg 199178, Russia \vskip-4pt {\small
nikolay.g.kuznetsov@gmail.com}

\vskip7mm

\parbox{136mm} {\noindent One isoperimetric inequality for the fundamental sloshing
eigenvalue is derived under the assumption that containers have vertical side walls
and either finite or infinite depth. It asserts that among all such containers,
whose free surfaces are convex, have two axes of symmetry and a given perimeter
length, this eigenvalue is maximized by infinitely deep ones provided the free
surface is either the square or the equilateral triangle. The proof is based on the
recent isoperimetric result obtained by A. Henrot, A. Lemenant and I.~Lucardesi for
the first nonzero eigenvalue of the two-dimensional Neumann Laplacian under the
perimeter constraint.

Another isoperimetric inequality for the fundamental eigenvalue, which describes
sloshing in containers with vertical walls, is a consequence of the classical
result due to G. Szeg\H o concerning the first nonzero eigenvalue of the free
membrane problem.}

\vskip4mm

{\centering \section{Introduction} }

\noindent 

Isoperimetric inequalities for eigenvalues of problems arising in mathematical
physics always attracted the interest of mathematicians being also important for
applications; see, for example, the survey article \cite{P} and the
monograph~\cite{B}. The best-known results in this area are as follows:

(i) The Faber--Krahn inequality (see \cite[p.~104]{B}), conjectured by Lord Rayleigh
in 1877 and proved in the 1920s, says that the ball minimizes the first eigenvalue
$\lambda_1$ of the Dirichlet Laplacian under the volume constraint on a domain.

(ii) The Szeg\H o--Weinberger inequality (see \cite[pp. 121 and 153]{B}) proved in
the 1950s says that the ball maximizes the first nonzero eigenvalue $\mu_1$ of the
Neumann Laplacian under the volume constraint on a domain.

(iii) The Weinstock inequality (see \cite[p.~463]{P}) also proved in the 1950s says
that the disk maximizes the first nonzero eigenvalue $\sigma_1$ of the Steklov
Laplacian under the perimeter constraint on a simply connected two-dimensional
domain.

Until recently, less was known about properties of the Neumann eigenvalue $\mu_1$
under the perimeter constraint, but some new results on this topic were obtained by
A. Henrot, A. Lemenant and I. Lucardesi \cite{HLL}. First, there is no maximizer of
$\mu_1$ in general, but it exists for two-dimensional, convex domains. Second, it is
conjectured that the inequality
\begin{equation}
P^2 (F) \, \mu_1 (F) \leq 16 \, \pi^2 \label{conj}
\end{equation}
is valid for every two-dimensional, convex domain $F$ ($P (F)$ denotes the perimeter
of $F$), and the equality is attained for squares and equilateral triangles;
moreover, this conjecture is proved for domains with two axes of symmetry, and it is
shown that squares and equilateral triangles are the only domains providing the
equality.
 
The latter assertion serves as the basis for obtaining an isoperimetric inequality
valid for the fundamental sloshing eigenvalue in the case of containers belonging to
a particular class; namely, those having vertical side walls, either finite or
infinite depth and the mentioned symmetry of a convex free surface.

The Szeg\H o inequality yields another isoperimetric inequality for the fundamental
sloshing eigenvalue in the case of containers with vertical walls.

Statement of the sloshing problem and some known facts about its eigenvalues are
collected in Sect.~2. New results are formulated and proved in Sect.~3.

\vspace{4mm}

{\centering \section{Sloshing problem} }

\subsection{Statement of the problem}

Here, the classical linear model of sloshing is used (see, for example, \cite{Lamb}
and \cite{HTW} for derivation of the governing equations), which as follows. Let an
inviscid, incompressible, heavy fluid, say water, occupies a three-dimensional
container $W \subset \RR^3_- = \{ (x,z) \in \RR^2 , \ y < 0 \}$. This Lipschitz
domain is bounded from above by a free surface, whose mean position while
oscillating (it is also its undisturbed position) coincides with a bounded,
two-dimensional domain $F$. Cartesian coordinates $(x,y,z)$ are chosen so that $F$
lies in the $(x,z)$-plane, whereas the $y$-axis is pointed in the direction opposite
to the gravitational force. The surface tension is neglected on $F$, and the water
motion is assumed to be irrotational and of small-amplitude.

These assumptions lead to the following boundary value problem for $\phi
(x,y,z)$\,---\,the velocity potential of the water motion with a time-harmonic
factor removed:
\begin{eqnarray}
&& \phi_{xx} + \phi_{yy} + \phi_{zz} = 0 \quad \mbox{in} \  W , \label{lap} \\ &&
\phi_y = \nu u \quad \mbox{on} \  F , \label{nu} \\ && \frac{\partial \phi}{\partial
n} = 0 \quad \mbox{on} \ \partial W \setminus \bar F . \label{nc}
\end{eqnarray}
Here, $\nu = \omega^2/g$ with $\omega$ being the radian frequency of the water
oscillations and $g$ standing for the acceleration due to gravity. It is convenient
to complement these relations by the orthogonality condition
\begin{equation}
\int_{F} \phi \, \D x \D z = 0 \, , \label{ort}
\end{equation}
thus excluding the zero eigenvalue of \eqref{lap}--\eqref{nc}.

In the case of unbounded $W$, a condition at infinity must be added for $\phi (p)$;
here and below $p = (x,y,z)$ for the sake of brevity. The following behaviour
\[ |\nabla \phi (p)| = o \left( |p|^{-2} \right) \quad \mbox{as} \ |p| = \big( x^2 
+ y^2 + z^2 \big)^{1/2} \to \infty 
\]
was proposed in \cite[p.~289]{HTW}. Being combined with \eqref{ort}, it guarantees
that
\begin{equation}
\int_{W} |\nabla \phi (p)|^2 \, \D x \D y \D z < \infty \, , \label{kener}
\end{equation}
that is, the sloshing motion has finite kinetic energy in an arbitrary unbounded $W$
(this is obvious for bounded $W$). An application of Green's formula shows that the
finiteness of potential energy
\begin{equation}
\int_F \phi^2 \, \D x \D z < \infty \, , \label{pener}
\end{equation}
is equivalent to \eqref{kener}.

\subsection{Variational principle}

It is well known that the sloshing problem in $W$ can be presented as a variational
problem and the corresponding Rayleigh quotient is as follows (see, for example,
\cite{M}):
\begin{equation}
R (\phi) = \frac{\int_W |\nabla \phi|^2 \,\D x \D y \D z}{\int_F \phi^2 \, \D x \D
z} . \label{Ray}
\end{equation}
To obtain the fundamental eigenvalue $\nu_1$, one has to minimize $R (\phi)$ over
the subspace of the Sobolev space $H^1 (W)$, consisting of functions that satisfy
the orthogonality condition \eqref{ort}. Moreover, applying minimax method to
\eqref{Ray}, one obtains the whole sequence $\{ \nu_n \}_1^\infty$ of sloshing
eigenvalues tending to $+ \infty$ as $n \to \infty$.

Furthermore, minimax method yields the following domain monotonicity principle (see
\cite{KLF}): {\it Let $W' \subset W$ be two bounded containers such that $F' = F$,
then the inequality}
\begin{equation}
\nu_n (W') \leq \nu_n (W) \ \ is \ valid \ for \ all \ n = 1, 2, \dots .
\label{monot}
\end{equation}
For $n = 1$ this result was obtained in \cite{M}.

{\centering \section{Isoperimetric properties of the fundamental eigenvalue} }

{\bf The case of volume constraint.} As early as 1965, B.~A. Troesch \cite{T}
demonstrated that a parabola furnishes the shape of an axially symmetric container
$W$, which lies strictly below $F$ of the given radius $r_0$, and has the fixed
volume $|W|$ and the largest fundamental sloshing eigenvalue. Let us summarize his
approach to illustrate its difference from that presented below.

Troesch seeks the $C^1$-function $h (r)$, given on $[0, r_0]$ and vanishing for $r =
r_0$, by minimizing $R (\phi)$ with $\phi = f (r) \cos (m \theta)$, where $(r,
\theta)$ are the polar coordinates in the $(x,z)$-plane. Moreover, it is required
that the expression
\[ \frac{\int_0^{r_0} h \, (f_r^2 + m^2 \, r^{-2} f^2) \, r \, \D r}{2 \pi 
\int_0^{r_0} h \, r \, \D r \int_0^{r_0} f^2 \, r \, \D r} \, , \quad m = 0,1,\dots
\, ,
\]
is stationary with respect to variation of $h$, because the latter solves of the
isoperimetric problem; here it is taken into account that $|W|$ is fixed. Combining
this fact and the Euler--Lagrange equation for $R$, one arrives at the problem:
\[ (r h)_r - m^2 h + \nu r^2 = 0 \, , \ \ r \in (0, r_0) ; \quad h_r (0) = 0 , \ \
h (r_0) =  0 \, .
\]
Since the lowest eigenvalue for the symmetry class $m=1$ is the lowest one for all
sloshing modes, this boundary value problem yields the following isoperimetric
bottom:
\[ h_{\rm iso} (r) = \nu_1 (r_0^2 - r^2) / 2 \, .
\]
Finally, it is shown that
\[ \nu_1 \leq \frac{8}{r_0^4} \int_0^{r_0} h (r) \, r \, \D r 
\]
is the corresponding isoperimetric inequality, in which equality occurs for $h_{\rm
iso} (r)$.

It is worth mentioning that a basin described by the function $h_{\rm iso} (r)$ with
an arbitrary parameter $h_0$ replacing $\nu_1$ appears also in the Lamb's treatment
of tidal waves; see \cite[p.~291]{Lamb}.

\vspace{2mm}

\noindent {\bf The case of perimeter constraint.} Let us consider containers that
have vertical side walls and simply connected free surfaces, whereas their depths
are either finite or infinite. Let us show that {\it the inequality
\begin{equation}
P (F) \, \nu_1 (W) \leq 4 \, \pi \label{iso}
\end{equation}
is valid for every such container $W$ with convex $F$ having two axes of symmetry.
Moreover, the equality in \eqref{iso} is attained only for infinitely deep $W$,
whose free surface $F$ is either square or equilateral triangle.}

\begin{proof}
Let $W$ have finite nonconstant depth, that is, $d = \max_{ p \in \partial W} |y| <
\infty$ and 
\[ W \subset W_d = F \times (-d, 0) , \quad W \neq W_d .
\]
Then the domain monotonicity principle \eqref{monot} implies that $\nu_1 (W) \leq
\nu_1 (W_d)$. Separation of variables yields that 
\[ \nu_1 (W_d) = \sqrt{\mu_1 (F)} \, \tanh (d \sqrt{\mu_1 (F)}) \, ,
\]
and so $\nu_1 (W_d) < \sqrt{\mu_1 (F)}$. Hence for arbitrary $W$ we have
\[ \nu_1 (W) \leq \sqrt{\mu_1 (F)} = \nu_1 (W_\infty) , \quad \mbox{where} \ W_\infty
= F \times (-\infty, 0) 
\]
and the last equality also follows by separation of variables. In view of
\eqref{conj} obtained in \cite{HLL} under the assumptions imposed on $F$, we have
that $P (F) \, \sqrt{\mu_1 (F)} \leq 4 \, \pi$ (which is \eqref{iso} for
$W_\infty$), and so inequality \eqref{iso} is valid for $W$ and $W_d$.

The assertion that equality in \eqref{iso} occurs only for infinitely deep $W$,
whose free surface~$F$ is either square or equilateral triangle follows from the
corresponding result about inequality \eqref{conj}; its proof is far from being
simple in \cite{HLL}.
\end{proof}


\noindent {\bf The case of area constraint.} There is another isoperimetric
inequality for containers having vertical side walls. It is a consequence of the
classical result due to Szeg\H o \cite{Sz} which says that \newline {\it the inequality
\begin{equation}
|F| \, \mu_1 (F) \leq |D| \, \mu_1 (D) = \pi \, \mu_1 (\DD) = \pi \, (j'_{1,1})^2
\approx 3.39 \, \pi  \label{Sz}
\end{equation}
is valid for any simply connected domain $F$ with equality holding only when $F$ is
a disk $D$.} Here, $\DD$ denotes the unit disk, $|\cdot|$ stands for the area of a
domain, $j'_{1,1}$ is the first positive zero of $J'_1$, and $J_1$ is the Bessel
function of order one. This allows us to formulate the following assertion. \newline
{\it The inequality
\begin{equation}
\sqrt{|F|} \, \nu_1 (W) \leq \sqrt{\pi \, \mu_1 (\DD)} = \sqrt\pi \, j'_{1,1} \approx
1,84 \, \sqrt\pi \label{isop}
\end{equation}
is valid for any container $W$ with vertical side walls and simply connected free
surface $F$. Moreover, the equality in \eqref{isop} is attained only for infinitely
deep $W$, whose free surface is $D$.}

\begin{proof}
As above, the domain monotonicity principle \eqref{monot} implies that if $W$ has
finite depth, then 
\[ \nu_1 (W) < \sqrt{\mu_1 (F)} = \nu_1 (W_\infty) \, .
\]
Therefore, for arbitrary $W$ we have
\begin{equation}
|F| \, [\nu_1 (W)]^2 \leq |F| \, \mu_1 (F) \, , \label{fin}
\end{equation}
where equality is valid only for infinitely deep $W$. Combining this and the Szeg\H
o inequality \eqref{Sz}, one arrives at \eqref{isop}.

Since equality in \eqref{fin} is valid only for infinitely deep $W$, whereas
equality in \eqref{Sz} is valid only for $D$, the equality in \eqref{isop} is
attained only for infinitely deep $W$, whose free surface is $D$.
\end{proof}

\renewcommand{\refname}{
\begin{center}{\Large\bf References}
\end{center}}
\makeatletter
\renewcommand{\@biblabel}[1]{#1.\hfill}
\makeatother

\end{document}